\documentclass[12pt,leqno]{article}
\usepackage{amsfonts,amssymb,fullpage}
\newtheorem{num}{\indent\hskip-4.5pt}[section]

\newcommand{\alku}{ \begin{num}\hskip-6pt \hskip5pt }
\newcommand{\loppu}{\end{num}}

\newcommand{\be}{\addtocounter{num}{1}\begin{equation}}
\newcommand{\ee}{\end{equation}}

\newcounter{minutes}
\divide\time by 60
\newcounter{hours}
\multiply\time by 60
\addtocounter{minutes}{-\time}

\begin{document}

\def\Rn{{\mathbb R}^n}
\def\Rnbar{\overline{{\mathbb R}^n}}
\def\Rbar{{\overline \mathbb R}}
\def\R{{\mathbb R}}
\def\Rp{{\mathbb R}^+}
\def\C{{\mathbb C}}
\def\X{{\mathbb X}}
\def\cross{\times}
\def\ple{\preceq}
\def\pge{\succeq}
\def\ch{{\cosh}}
\def\sh{{\sinh}}
\def\arsh{{\rm arsh }}
\def\arch{{\rm arch }}
\def\lem{{\rm {\bf Lemma.\hskip 0.5truecm}}}
\def\pro{{\bf Proposition.\hskip 0.5truecm }}
\def\cor{{\bf Corollary.\hskip 0.5truecm}}
\def\thm{{\bf Theorem.\hskip 0.5truecm}}
\def\rem{{\bf Remark.\hskip 0.5truecm}}
\def\define{{\bf Definition.\hskip 0.5truecm}}
\def\ack{{\bf Acknowledgement.\hskip 0.5truecm}}
\def\lqq{\lq\lq}
\def\rqq{\rq\rq\,}
\def\rqqs{\rq\rq\hskip 0.15truecm}
\def\ineq{\not=}
\def\hbar{\vert}
\def\l{\ell}
\def\inv{ {\rm inv} } 
\def\abrv{.\hskip 0.1truecm}
\def\ident{\equiv}
\def\card{{\rm card\, }}
\def\GM{{GM(B^n)}}
\def\qed{$\square$}

\pagestyle{myheadings}
\title{Inequalities of relative weighted metrics}
\author{Peter A. H\"ast\"o,\\ 
{\small {Department of Mathematics, University of Helsinki, 
P.O. Box 4, 00014, Helsinki, Finland,}}\\ 
{\small{E-mail: peter.hasto@helsinki.fi.}} }
\date{\small{3rd of January, 2002. Bullv9.tex}}
\maketitle

{
\renewcommand{\thefootnote}{}
\tiny \footnote{\centerline{Mathematics Subject Classification (2000): Primary 51M16.}}
}
\begin{abstract}
In this paper we present inequalities between two generalizations of 
the hyperbolic metric and the $j_G$ metric. We also prove inequalities 
between generalized versions of the $j_G$ metric and Seittenranta's metric.
\end{abstract}



\section{Introduction}
\label{intro}

This paper contains various inequalities between metrics defined in 
subdomains $G$ of the M\"obius space $\Rnbar := \Rn \cup \{ \infty \}$, 
$n\ge 2$. In what follows all topological operations are with respect to $\Rnbar$ 
(see Section 2, for further reference e.g\abrv \cite{Vu}). 
We will always denote by $G\subset\Rnbar$ a 
domain (i.e\abrv open and connected set) with at least two boundary points
and by $x$ and $y$ points in $G$ similarly for $G'$, $x'$ and $y'$. 

This section contains the definitions of the metrics studied as well as the 
statement of the main results. The main results are two blocks of inequalities, 
Theorems \ref{rhoIneq} and \ref{delta_pIneq} which concern two different 
generalizations of the hyperbolic metric. Section 2 describes the notation 
used in this paper, which conforms to that used in \cite{Vu}. 
The two main theorems are proved in Sections 3 and 4, respectively. 

Our first result is a comparison between the  
generalized hyperbolic metric which was introduced in {\cite[(3.28)]{Vu}}, and proven 
to be a metric in domains with at least two boundary points 
in {\cite{Ha2}} 
with the generalized hyperbolic metric introduced by 
Pasi Seittenranta in \cite[Definition 1.1]{Se} and the well-known $j_G$ metric defined 
for $G\subset \Rn$ by
$$ j_G(x,y) := \log\left( 1 + {|x-y| \over \min\{ d(x), d(y)\}} \right).$$
For simplicity, the generalized hyperbolic metric from \cite{Vu} will 
be called {\it the generalized hyperbolic metric} or the {\it $\rho_G$ metric}, 
whereas that from \cite{Se} will be called {\it Seittenranta's metric} or the 
{\it $\delta_G$ metric}. For domains $G$ with at least two boundary points, 
the generalized hyperbolic metric is defined by 
\be\label{rhoDefIntro} 
\rho_G(x,y):= \sup_{a,b\in \partial G} \arch\{1 + |a,x,b,y||a,y,b,x|/2 \}\ee 
and Seittenranta's metric is defined by 
$$ \delta_G(x,y) := \sup_{a,b \in \partial G} \log \{1 + \vert a,x,b,y \vert \}, $$
where $\vert a,x,b,y \vert$ denotes the cross-ratio, see (\ref{crDef}). 

We cite some basic desirable properties of $\rho_G$ from \cite{Vu}, as this may help 
motivate studying this metric. Note that
$\delta_G$ also has all of these properties except the third which 
is replaced by $\delta_G(x,y)\ge \exp\{(q(\partial G)q(x,y)) \} -1$.
(Theorem 3.1 and Remark 3.2(2) of \cite{Se}, \cite[8.38(3)]{Vu}) 

\alku\lem\label{basicsLem} {\rm (\cite[3.25 \& 3.26]{Vu})} 
\begin{itemize}
\item[{\rm(i)}] $\rho_G$ is M\"obius invariant (see \cite[p. 32]{Be}). 
\item[{\rm(ii)}]$\rho_G$ is monotone in $G$, that is, if $G\subseteq G'$ then 
$\rho_{G'}(x,y)\le \rho_G(x,y)$ for all $x,y\in G$. 
\item[{\rm(iii)}] $\rho_G(x,y)\ge \cosh\{(q(\partial G)q(x,y))^2\} -1 $. 
\item[{\rm(iv)}] For $G=B^n$ and $G=H^n$ (the upper half-space), $\rho_G$ equals the hyperbolic metric. 
\end{itemize}
\loppu 

In this paper we prove the following inequalities of $\rho_G$:

\alku\thm\label{rhoIneq}
Let $G$ be a domain with $\card \partial G \ge 2$. Then
\begin{itemize}
\item[\rm{(i)}] $\delta_G \le \rho_G \le \frac{\arch\, 3}{\log 3} \delta_G$.
\end{itemize}
Assume additionally that $G\subset \Rn$. Then
\begin{itemize}
\item[\rm{(ii)}] $ j_G \le \rho_G \le {\arch\, 3\over \log 2} j_G$.
\end{itemize}
Both inequalities in (i) and the former inequality in (ii) are sharp.
\loppu 

\alku\rem{\rm Note that the term \lqq sharp\rqq when applied to an inequality 
means that the constant cannot be improved, i.e\abrv there exists 
points $x_i,y_i\in G$, $i=1,2, ...$, such that 
$$ \lim_{i\to\infty} d_1(x_i,y_i)/d_2(x_i,y_i) = c,$$
for the inequality $d_1\le c d_2$. 
}\loppu

It was shown in \cite[Corollary 6.1]{Ha1} that $\delta_G$ can be 
embedded in the following family of metrics ($0<p<\infty$):
$$ \delta_G^p(x,y):= \sup_{a,b\in \partial G} \log \{1 + (\vert x,a,y,b\vert^p+ 
\vert x,b,y,a\vert^p)^{1/p} \},\ 
\delta_G^{\infty}(x,y):= \lim_{p\to\infty} \delta_G^p(x,y).$$
With this notation $\delta_G^\infty=\delta_G$, Seittenranta's metric. 
It likewise follows directly from Remark 6.1 and Corollary 6.1 in 
\cite{Ha1} that for $G\subset \Rn$, $j_G$ can be embedded in the family 
$$ j_G^p(x,y):= \sup_{a\in \partial G} \log \left( 1 + \left( 
{|x-y|^p \over |x-a|^p} +{|x-y|^p \over |y-a|^p}\right)^{1/p} \right),\ 
j_G^{\infty}(x,y):= \lim_{p\to\infty} j_G^p(x,y).$$
where $0<p<\infty$. Here then $j_G^\infty=j_G$, the classical $j_G$ metric.
We note that if we fix $b=\infty$ in the expression for 
$\delta_G^p$ then we get the expression for $j_G^p$.

In this paper we prove the following inequalities of the generalized 
$j_G$ and $\delta_G$ metrics. 

\alku\thm\label{delta_pIneq} Let $G$ be a domain with at least two boundary 
points. If $0 < q\le p\le \infty$ then
\begin{itemize}
\item[\rm{(i)}] $\delta_G^p\le \delta_G^q \le 2^{1/q-1/p} \delta_G^p$. 
\end{itemize}
If additionally $G\subset \Rn$ then 
\begin
{itemize}
\item[\rm{(ii)}] $j_G^p\le j_G^q \le 2^{1/q-1/p} j_G^p$.
\end{itemize}
If $p\in[1,\infty]$ and $G\subset \Rn$ then 
\begin{itemize}
\item[\rm{(iii)}] $j_G^p\le \delta_G^p\le 2j_G^p$.
\end{itemize}
All the inequalities are sharp.
\loppu

Note that inequality (iii) of the previous theorem 
is a generalization of \cite[Theorem 3.4]{Se}.


\section{Notation}
\label{definitions}

The notation adopted here mostly corresponds to that of \cite[Chapter 2]{Vu}, the same 
material is also presented in e.g\abrv \cite[Chapter 3]{Be}. 
We denote by $\{e_1, e_2, ..., e_n\}$ the 
standard basis of $\Rn$ and by $n$ the dimension of the Euclidean space
under consideration and assume that $n\ge 2$. For $x\in\Rn$ we denote by 
$x_i$ the $i^{\rm th}$ coordinate of $x$. 
We will identify $\R$ with the subspace $\R e_1$ of 
$\Rn$. Hence if $x$ is a real number then the expression \lqq the point $x$\rqq means the 
point $x e_1$ etc. We will use the notation $\Rbar := \R \cup \{ \infty, -\infty\}$ and 
$\Rnbar := \Rn \cup \{ \infty\}$ for the two and one point compactifications of 
$\R$ and $\Rn$, respectively. 

By $\partial G$ we will denote the boundary and by $G^c$ the complement of $G$ 
with respect to $\Rnbar$. The following notation will be used for balls, spheres and 
the upper half-space:
$$ B^n(x, r):=\{y \in \Rn \colon |x-y| <r \},\quad S^{n-1}(x, r):=\partial B^n(x,r),
\quad H^n:=\{x\in\Rn\colon x_n>0\}. $$

We define the spherical metric $q$ in $\Rnbar$ 
by means of the canonical projection onto the Riemann sphere, hence
\be\label{qmetric} q(x,y):={ \vert   x-y \vert   \over  \sqrt{1+\vert  x\vert^2 
}\sqrt{1+\vert y\vert^2 }},\ q(x,\infty):= {1 \over \sqrt{1+\vert  x\vert^2 }}. \ee
We will consider $\Rnbar$ as the metric space $(\Rnbar, q)$, hence its balls are the 
(open) balls of $\Rn$ and complements of closed balls of $\Rn$ as well as half-spaces.
The cross-ratio $\vert a,b,c,d\vert$ is defined by
\be\label{crDef} |a,b,c,d| :={q(a,c)q(b,d)\over q(a,b) q(c,d)}\ee
for $ a,b,c,d \in\Rnbar$, $a\ineq b$ and $c\ineq d$. If $ a,b,c,d \in\Rn$ then 
the cross-ratio can be expressed in terms of Euclidean distances as
$$ |a,b,c,d| :={|a-c||b-d|\over |a-b||c-d|}.$$ 
A mapping $f\colon \Rnbar \to \Rnbar$ is a M\"obius mapping if 
$$ |f(a),f(b),f(c),f(d)| = |a,b,c,d| $$
for every quadruple $a,b,c,d\in \Rnbar$ with $a\ineq b$ and $c\ineq d$ 
(\cite[p. 32]{Be}). 

We denote the inverses of the hyperbolic sine and cosine by
$$ \arsh(x)=\log(x+\sqrt{x^2+1}),\ x\ge 0, $$
and 
$$ \arch(x)=\log(x+\sqrt{x^2-1}),\ x\ge 1, $$
respectively. Note that $\sh(\arch(x))=\sqrt{x^2-1}$.


\section{The proof of Theorem 1.3}
\label{other_ineqs}

{\bf Proof of Theorem 1.3(i). \hskip 0.5truecm}
We start by proving the first of the inequalities. 
Fix the points $x,y\in G$ and $a,b\in\partial G$ such that 
$\delta_G(x,y)=\log\{1+|a,x,b,y|\}$. The points $a$ and $b$ can be chosen, since 
$\partial G$ is a compact set in $\Rnbar$. Then it suffices to prove the first 
inequality in 
\be\label{rIe1}  \log\{1+|a,x,b,y|\} \le \arch\{ 1 + 
\vert a,x,b,y \vert \vert a,y,b,x \vert/2 \} \le \rho_G(x,y),\ee
since the second follows directly from the definition of $\rho_G$. 
Moreover, since both $\delta_G$ and $\rho_G$ are M\"obius invariant 
we may assume that $a=\infty$ and $b=0$. Denote 
$s:= \vert x-y\vert /\sqrt{\vert x\vert \vert y\vert}$ and 
$k:= \sqrt{\vert x\vert / \vert y\vert }$ and assume that 
$\vert x\vert \ge \vert y\vert$. Then (\ref{rIe1}) becomes
$$ \log \{1+ k s\} \le 
\log \{1+s^2/2+\sqrt{s^4/4+s^2}\} $$
which reduces to $k- s/2 \le \sqrt{s^2/4+1}$. Squaring this, we see
that the inequality holds, since $s \ge k - 1/k $ by the definitions of $k$ and $s$ using 
the Euclidean triangle inequality. We see that there is equality for 
$G=\Rn\setminus\{0\}$, $x=e_1$ and $y=r e_1$, $r\in \R$.

In proving the second inequality it again suffices to assume $a=\infty$ and $b=0$.
Let $s$ and $k$ be as before and set $c:= \arch\{3\}/ \log\{3\}$. The second 
inequality is equivalent to
\be\label{iTe2} c \log \{ 1+ ks\} - \log \{1+s^2/2+\sqrt{s^4/4+s^2}\}\ge 0. \ee
The derivative with respect to $s$ of the left hand side of the above inequality equals
$$ {c\over s+ 1/k } - {1 \over \sqrt{s^2/4+1}} =
{1\over s+ 1/k } \left( c - {s+ 1/k \over \sqrt{s^2/4+1}} \right).$$
Since the term in the parenthesis is decreasing in $s$, the derivative has at most
one zero, which is a maximum. Therefore we need only check that (\ref{iTe2}) holds at the 
end-points, $s=0$ and $s=k+1/k$, which correspond to $|x-y|=0$ and $|x-y|=|x|+|y|$, 
respectively. For $s=0$ the inequality (\ref{iTe2}) obviously holds. 
In the case $s=k+1/k$, since $k = s/2 + \sqrt{s^2-1}/2$, we need to show that
$$ c\log \{1 + s^2/2 + \sqrt{s^4/4-s^2} \} \ge \log \{1 + s^2/2 + \sqrt{s^4/4+s^2}\}.$$
Clearly equality holds for $s=2$. The claim then follows when we show that the 
left hand side has greater derivative than the right hand side for $s\ge 2$. 
Let us change variable, $t=s^2$, and differentiate with respect to $t$:
$$ c {t^2 + t\sqrt{t^2-4} -1 \over \sqrt{t^2-4}(2t^2+1)} \ge {1\over \sqrt{t^2+4}}.$$
Since $c\ge 1$, we may drop it. Multiplying by $\sqrt{t^2-4}(2t^2+1)\sqrt{t^2+4}$ 
and squaring gives, after rearranging and dividing by $2$, 
$$ t(t^2-1)(t^2+4)\sqrt{t^2-4} \ge t^6 - 8t^4 + 4t^2 - 4.$$
To see that this holds, observe the following chain of inequalities (note 
that $t\ge 4$): 
$$ t(t^2-1)(t^2+4)\sqrt{t^2-4} \ge t^5\sqrt{t^2-4} \ge t^6 - 7t^4 \ge 
t^6 - 8t^4 + 4t^2 - 4. $$

For the sharpness of this inequality we choose $G=\Rn\setminus\{0\}$, $x=1$ and $y=-1$. 
Then there is equality in the inequality, and hence the constant cannot be improved.
\qed

\bigskip

{\bf Proof of Theorem 1.3(ii). \hskip 0.5truecm} 
The first inequality follows from Theorem \ref{rhoIneq}(i) ($\delta_G\le\rho_G$) 
and \cite[Theorem 3.4]{Se} ($j_G\le\delta_G$). Its sharpness follows by 
letting $G=H^n$, $x= se_n$ and $y=re_n$ (see \cite[Remark 3.5]{Se}).

We turn to the second inequality. The metric $j_G$ as it is normally defined is 
not M\"obius invariant, and indeed $\infty$ is a special point in the sense that it 
may not belong to the domain $G$ in which the metric is defined. We may, however, 
think of $j_G$ as the member $j_{G,\infty}$ of the following family:
\be\label{jMobius} j_{G,b}(x,y):= \sup_{a\in \partial G} \log \{ 1 + 
\max\{ \vert x,a,y,b\vert, \vert x,b, y, a\vert\} \},\ee
where $G\subset \Rnbar$ is a domain not containing $b$ with at least two boundary points. 
Since $j_{G,b}$ is defined in terms of cross ratios, it is clear that it is 
M\"obius invariant. Hence we may apply an auxiliary M\"obius transform to 
both sides of the inequality $\rho_G\le j_{G,\infty}$, as long as we keep track of 
where $\infty$ is mapped and use the appropriate $j_{G,b}$.

As before we may then assume that the boundary points $a$ and $b$ 
occurring in the definition of $\rho_G$ equal $0$ and $\infty$. 
We need to prove
\be\label{T2E1} \arch\left(1+ {\vert x'-y'\vert^2 \over 2\vert x'\vert 
\vert y'\vert}\right) \le {\arch\, 3\over \log 2} 
\sup_{a\in \partial G'} \log \{ 1 + \max\{\vert \infty,a,y',b\vert, 
\vert \infty, b, y', a\vert\}\} \ee
$$ \le {\arch\, 3\over \log 2} \sup_{a\in \partial G'} \log \left( 1 + 
{|x'-y'| |b-a| \over \min\{|x'-b||y'-a|, |x'-a||y'-b|\}}\right), $$
where the supremum is over the boundary point $a$ only; $b$ is some fixed 
point in the complement of $G$. 
If $b=0$ or $b=\infty$, we may proceed exactly as in the proof of Theorem 
\ref{rhoIneq} (i) and arrive at the better constant $\arch\{3\}/\log\{3\}$. Assume 
then that $b\not\in\{0,\infty\}$. 
We may then assume without loss of generality that $b=1$ (recall that $1$ denotes 
the point $e_1$) by scaling and rotating. 
Since both sides are M\"obius invariant, we may assume that $|x'||y'|\le 1$ 
by performing an inversion in the unit sphere, since this leaves $b$ fixed. We then 
forget about the original $x$ and $y$ and denote $x'$ by $x$ and $y'$ by 
$y$, to simplify the notation. 

We may restrict the supremum from $\partial G$ to $\{0,\infty\}\subseteq \partial G$ in 
the right hand side of (\ref{T2E1}), since this only makes the supremum smaller. 
Moreover, we can move the supremum to the denominator of the fraction inside the 
logarithm, changing it to minimum, since it is taken over finitely many terms. 
We then see that it suffices to prove that  
$$ \sup_{a\in \partial G} \log \left( 1 + 
{|x-y| |a-1| \over \min\{|x-1||y-a|, |x-a||y-1|\}}\right) \ge $$
$$ \ge \log \left( 1 + {|x-y| \over \min\{|x-1||y|, |x||y-1|, |y-1|, |x-1|\}}\right).$$
Let us estimate 
$|y-1| \le |y|+1$ and $|x-1|\le |x|+1$. Then 
$$ \min\{|x-1||y|, |x||y-1|, |y-1|, |x-1|\} \le $$
$$ \le \min\{ \min\{1,|y|\} (1+|x|) , \min\{1,|x|\} (1+|y|)\} 
= \min\{|x|,|y|\} + \min\{ 1, |x||y|\}.$$
Recall that we assumed that $|x||y|\le 1$. By symmetry, we may assume that $|x|\le |y|$. 
Then we need to prove
\be\label{T2E2} 
\arch\left(1+ {\vert x-y\vert^2 \over 2\vert x\vert \vert y\vert}\right) 
\le {\arch\, 3\over \log 2} \log \left( 1 + {|x-y| \over |x| + |x||y|} \right).\ee

Denote $S_c := \{z\in\Rn\colon |z-y|=c|z|\}$. For fixed $y$ and $c>0$ 
consider how the inequality (\ref{T2E2}) varies as $x$ varies over $S_c$:
\be\label{S_cIneq} \arch\left(1+ {c |x-y| \over 2 |y| }\right) \le 
{\arch\, 3\over \log 2} \log \left( 1 + {c\over 1 + |y|} \right).\ee
We see that the right hand side does not depend on $|x-y|$, which 
means that it suffices to consider points $x\in S_c$ which maximize this 
distance, since this yields the \lqq hardest\rqq inequality. 

Observe that for all $c>0$ the sphere $S_c$ intersects the
segment $[0,y]$ and $S_c$ encloses $y$ if and only if $c \in
(0,1)$ and $0$ if and only if $c>1$. Note also that $S_c$ is a subset
of  $B^n(|y|)$ if and only if $c>2$. Since we need only consider points $x$ that 
satisfy $|x| \le |y|$, we see that for $c \in(0,2]$, the distance $|x-y|$
is maximized by some $x$ satisfying $|x|=|y|.$ If $c>2$ $|x-y|$ is maximized by the
choice $x=-y/(c-1)$.

Let $\lambda:= \sqrt{|x||y|}\le 1$. If $\lambda <1$ then we can consider the points 
$x':=x/\lambda $ and $y':=y/\lambda $. The left hand side of (\ref{T2E2}) is the 
same for the points $x$ and $y$ as for the points $x'$ and $y'$, however the 
right hand side is smaller for the latter points. Hence we see that it suffices 
to prove (\ref{T2E2}) for points $x$ and $y$ with $|x||y|=1$.

Combining the conclusions of the previous two paragraphs, we see that 
if $c\le 2$ we need to consider only the case $|x|=|y|=1$, i.e. 
$$ \arch\{1+ s^2 /2\} \le {\arch\, 3\over \log 2} \log\{ 1 + s/2\}, $$
where we have denoted $s:=|x-y|$. For $s=0$ there is equality in the 
inequality, and since the left hand side has lesser derivative than the 
right hand side we are done with this case.

In the case $c<2$ we need to consider points $x$ and $y$ with $|x||y|=1$ 
such that $x$, $0$ and $y$ lie on some line in this order. Hence we need to show that 
$$ \arch\left(1+ {(|x|+|y|)^2 \over 2|x| |y|}\right) \le 
{\arch\, 3\over \log 2} \log \left( 1 + {|x|+|y| \over |x| + |x||y|} \right).$$
Let us write $t:=|y|=1/|x|\ge 1$. The previous inequality becomes 
$$ \arch\{1+ (t+1/t)^2/2 \} \le {\arch\, 3\over \log 2} \log \{1 + (t^2+1)/(t+1)\}.$$
For $t=1$ there is clearly equality in this inequality. We 
show that the right hand side has larger derivative than the left hand side 
for all $t>1$, which is equivalent to 
$$ {2\over t} {t^2-1\over \sqrt{t^4+6t^2+1}} \le {\arch\, 3\over \log 2} 
{t^2+2t-1 \over (t+1)(t^2+t+2)}.$$
We use the estimate $\arch\, 3 /(2 \log 2) \ge 5/4$ and multiply both sides by the denominators:
$$ 4(t^2-1) (t+1)(t^2+t+2) \le 5 (t^2+2t-1)t\sqrt{t^4+6t^2+1}.$$
We then use the estimates $t^2+2t-1 \ge t(t+1)$  and $\sqrt{t^4+6t^2+1} \ge t^2+1$
and cancel the term $t+1$ from both sides: 
$$ 4(t^4+t^3+t^2-t-2)\le 5(t^4+t^2).$$ 
With the substitution $u:=t+1$ this is equivalent to $u^4-5u^2-2u+10 \ge 0$. 
Since $2u\le u^2+1$ we have $u^4-5u^2-2u+10 \ge u^4 -6u^2+9\ge (u^2-3)^2\ge 0$.
\qed


\section{Proof of inequality 1.5}
\label{rel_ineqs}

{\bf Proof of Theorem 1.5 (i) and (ii). \hskip 0.5truecm}
It suffices to prove each of the claims for some fixed boundary point(s), 
since we may choose it (them) to correspond to the point(s) where the supremum is 
attained in the quantity whose upper bound we want to establish. Hence it suffices to 
prove the real-number inequality
$$ \log(1+(x^p+y^p)^{1/p})\le \log(1+(x^q+ y^q)^{1/q})\le 
2^{1/p-1/q} \log(1+(x^p+y^p)^{1/p})$$
in order to prove both of the claims. 
Since $(x^p+y^p)^{1/p} \le (x^q+y^q)^{1/q}$ the first 
inequality is clear. Let us denote $s:=2^{1/q-1/p}\ge 1$. Then 
$\log(1+xs)\le s\log(1+x)$ for $x\ge 0$ by the Bernoulli inequality. Hence 
it suffices to prove the first inequality in
$$ \log(1+(x^q+ y^q)^{1/q})\le \log(1+s(x^p+y^p)^{1/p}) \le s\log(1+(x^p+y^p)^{1/p}).$$
However, this is immediately clear, since $(x^q+y^q)^{1/q}\le s (x^p+y^p)^{1/p}$ 
by the power-mean inequality. 

We still need to show that the inequalities are sharp: Let $G:=\Rn\setminus \{0\}$. Then 
$$ \delta_G^p(x,y)=j_G^p(x,y) = \log \left( 1 + \left( 
{|x-y|^p \over |x|^p} +{|x-y|^p \over |y|^p}\right)^{1/p} \right).$$
Fix $y$ and let $x\to \infty$. Then 
$$ \lim_{ x\to \infty} {j_G^p(x,y)\over \log |x|} \to \log 2$$ 
irrespective of the value of $p$, which 
shows that the first inequalities are sharp. If $|x|=|y|$ then 
$$ \delta_G^p(x,y) = \log (1+ 2^{1/p} |x-y|/|x|).$$
As $x\to y$ we see that the second inequalities are also sharp. \qed

\alku\rem{\rm 
Note that since we are using a point-wise estimate, we need not consider the 
cases $p=\infty$ and $q=\infty$ separately, since $j_G^p$ and $\delta_G^p$ are 
both extended to $p=\infty$ continuously.
}\loppu

{\bf Proof of Theorem 1.5 (iii). \hskip 0.5truecm}
The first inequality follows since $b$ can equal $\infty$ in 
the definition of $\delta_G$, in which case $\delta_G^p=j_G^p$. In the 
domain $\Rn\setminus\{0\}$ we have $\delta_G^p(x,y)=j_G^p(x,y)$ for every pair of 
points $x,y\in G$, hence the inequality is sharp. 
It remains to consider the second inequality.

Fix $x$ and $y$ in $G$ and the boundary points $a$ and $b$ for which the 
supremum is attained. We may assume without loss of generality that 
$|x-y|=1$. Then
$$ \delta_G^p(x,y)= \log \{1+(\vert x,a,y,b\vert^p+ \vert x,b,y,a\vert^p)^{1/p}\}\le $$
$$ \le \log \{ 1+ ((s+t+st)^p + (u+v+uv)^p)^{1/p}\}, $$
where we have denoted 
$$ s:= {1\over |x-a|},\; t:= {1\over |y-b|}\; u:= {1\over |x-b|}\; v:= {1\over |y-a|}, $$
and used the estimates
$$ |a-b|\le |a-x|+|x-y|+|y-b|{\rm\ and\ }|a-b|\le |a-y|+|y-x|+|x-b|. $$ 
in $|x,a,y,b|$ and $|x,b,y,a|$, respectively. Now 
$$ j_G^p(x,y) \ge \sup_{w\in \{a,b\}} \log\{1+ (|x-w|^{-p} +|y-w|^{-p})^{1/p}\} = \log\{1+\max\{s^p+v^p, t^p+u^p\}^{1/p} \}. $$
By symmetry we may assume that the $t^p+u^p \le s^p+v^p$. If we apply the exponential function to both sides of the inequality
\begin{eqnarray*} 
\delta_G^p(x,y) & \le & \log \{ 1+ ((s+t+st)^p + (u+v+uv)^p)^{1/p}\} \\
& \le & 2\log\{1+\max\{s^p+v^p, t^p+u^p\}^{1/p} \} \le 2 j_G^p(x,y),
\end{eqnarray*}
we see that it suffices to show that 
\be\label{jd2jE1} 1+ ((s+t+st)^p + (u+v+uv)^p)^{1/p} \le 
(1+ (s^p+v^p)^{1/p})^2. \ee
We see that the left hand side can be increased by increasing $t$ while keeping the 
right hand side constant if $t^p+u^p< s^p+v^p$. Hence we may assume that 
$t^p+u^p=s^p+v^p=:\alpha^p$.

We will show that (\ref{jd2jE1}) holds for every quadruple $s,t,u,v\in\Rp$ 
with $t^p+u^p=s^p+v^p$ for $p\ge 1$. For fixed $s$, $t$, $u$ and 
$v$ let us consider how the inequality varies under the transformation 
$x\mapsto wx$, $y\mapsto wy$, $u\mapsto wu$ and $v\mapsto wv$.
Then the equation (\ref{jd2jE1}) becomes, after we divide it by the common factor $w$,
$$ f(w):=2\alpha + w\alpha^2 - ((s+t+stw)^p + (u+v+uvw)^p)^{1/p}\ge 0.$$
We will show that $f$ increases in $w$. The derivative $f'(w)$ equals 
$$ \alpha^2 - \{ (s+t+stw)^p + (u+v+uvw)^p\}^{1/p-1} 
\{(s+t+stw)^{p-1}st + (u+v+uvw)^{p-1}uv\} = $$ 
$$ = \alpha^2 - \{(1+\zeta^p)^{1/p-1} st + (1+\zeta^{-p})^{1/p-1} uv\}=:
\alpha^2-g(\zeta),$$
where $\zeta := (u+v+uvw)/(s+t+stw)$. We will now consider how 
$$ g(\zeta)=(1+\zeta^p)^{1/p-1} st + (1+\zeta^{-p})^{1/p-1} uv $$ 
varies with $\zeta$. The derivative $g'(\zeta)$ equals 
$$ -(p-1)( (1+\zeta^p)^{1/p-2} \zeta^{p-1} st-(1+\zeta^{-p})^{1/p-2}\zeta ^{-p-1} uv)= $$
$$ = -(p-1) (1+\zeta^p)^{1/p-2} \zeta ^{p-2}(st\zeta - uv).$$
We see that $g$ has a maximum at $\zeta=uv/(st)$ for $p>1$. Hence 
\begin{eqnarray*} 
{df\over dw} & \ge & \alpha^2 -g\left({uv\over st}\right)= 
\alpha^2 - \left( \left( {(st)^p+(uv)^p\over (st)^p }\right)^{1/p-1} st + 
\left( {(st)^p+(uv)^p\over (uv)^p} \right)^{1/p-1}  uv\right)\\
& = & \alpha^2- ((st)^p+(uv)^p)^{1/p} = 
(t^p+u^p)^{1/p}(s^p+v^p)^{1/p}-((st)^p+(uv)^p)^{1/p}\ge 0.
\end{eqnarray*}
Now since $f$ is increasing in $w$, it suffices to show that $f(0)\ge 0$ 
in order to obtain $f(w)\ge 0$, which is equivalent with (\ref{jd2jE1}). 
In other words we must show that 
$$ 2(s^p+v^p)^{1/p}  - ((s+t)^p + (u+v)^p)^{1/p}\ge 0.$$
Recall that $t^p+u^p=s^p+v^p=:\alpha^{1/p}$ and denote additionally $\beta:=s+t$. The 
previous inequality becomes
$$ 2\alpha - \{\beta^p + [(\alpha^p-s^p)^{1/p} + 
(\alpha^p-(\beta-s)^p)^{1/p}]^p\}^{1/p}.$$
For fixed $\alpha$ and $\beta$, $(\alpha^p-s^p)^{1/p} + (\alpha^p-(\beta -s)^p)^{1/p} \le 2(\alpha^p-(\beta /2)^p)^{1/p}$ and so it suffices to show that 
$2\alpha - (\beta^p + 2^p(\alpha^p-(\beta/2)^p))^{1/p}\ge 0$, which is obvious. 

We still have to show that the inequality is sharp. Consider then the domain
$G=\Rn\setminus\{-e_1,e_1\}$ and the point $\epsilon e_2$ and $-\epsilon e_2$. We have
$$ \delta_G^p(\epsilon e_2, -\epsilon e_2) = 
\log \left(1+{2^{1/p+2}\epsilon \over \sqrt{1+\epsilon^2}}\right)$$
and 
$$ j_G^p(\epsilon e_2, -\epsilon e_2) = 
\log \left(1+{2^{1/p+1}\epsilon \over \sqrt{1+\epsilon^2}}\right).$$
It is then clear that 
$$ \lim_{\epsilon\to 0}{\delta_G^p(\epsilon e_2, -\epsilon e_2) \over 
j_G^p(\epsilon e_2, -\epsilon e_2)} = 2.$$

\alku\rem {\rm It is not immediately clear whether the inequality from 
Theorem \ref{delta_pIneq} (iii) holds for $0<p<1$ as well. It is clear that 
(\ref{jd2jE1}) does not hold in this case for arbitrary $s,t,u,v\in\Rp$, however, 
these variables are not really arbitrary but rather related 
by various triangle inequalities. 
}\loppu

\ack {\rm 
I would like to thank Matti Vuorinen for numerous comments and suggestions 
during the life-span of this manuscript as well as Pentti J\"arvi for his 
comments on a previous version of the manuscript.
}

\end{document}